\documentclass{amsart}
\usepackage{amsfonts,amsmath}
\usepackage{amssymb,color}

\begin{document}

\title[Self-Dual Maps]{Self-Dual Maps and Symmetric Bistochastic Matrices}
\author[Corey O'Meara]{Corey O'Meara}
\address{Department of Chemistry, Technical University of Munich, 85747 Garching, Germany}
\email{comeara@tum.de}
\author[Rajesh Pereira]{Rajesh Pereira}
\address{Department of Mathematics \& Statistics, University of Guelph, Guelph, ON, Canada N1G 2W1}
\email{pereirar@uoguelph.ca}


\begin{abstract}

The $H$-unistochastic matrices are a special class of symmetric bistochastic matrices obtained by taking the square of the absolute value of each entry of a Hermitian unitary matrix.  We examine the geometric relationship of the convex hull of the $n$ by $n$ $H$-unistochastic matrices relative to the larger convex set of $n$ by $n$ symmetric bistochastic matrices. We show that any line segment in the convex set of the $n$ by $n$ symmetric bistochastic matrices which passes through the centroid of this convex set must spend at least two-thirds of its length in the convex hull of the $n$ by $n$ $H$-unistochastic matrices when $n$ is either three or four and we prove a partial result for higher $n$.  A class of completely positive linear maps called the self-dual doubly stochastic maps is useful for studying this problem.  Some results on self-dual doubly stochastic maps are given including a self-dual version of the Laudau-Streater theorem.

\end{abstract}
\maketitle

\renewcommand{\vec}[1]{\mbox{\boldmath$#1$}}  
\newcommand{\e}{\mathrm{e}}                   
\newcommand{\dif}{\mathrm{d}}                 
\renewcommand{\span}[1]{\mathrm{span}\left\{#1\right\}}  
\newcommand{\N}{\mathbb{N}}                   
\newcommand{\Z}{\mathbb{Z}}                   
\newcommand{\Q}{\mathbb{Q}}                   
\newcommand{\R}{\mathbb{R}}                   
\newcommand{\C}{\mathbb{C}}                   
\renewcommand{\L}{\mathcal{L}^2(\R)}          
\renewcommand{\l}{\ell^2(\Z)}                
\renewcommand{\bar}[1]{\overline{#1}}         

\newtheorem{theorem}{Theorem}
\newtheorem{corollary}[theorem]{Corollary}
\newtheorem{proposition}{Proposition}
\newtheorem{definition}{Definition}

\newtheorem{notation}{Notation}

\newtheorem{lemma}{Lemma}

\newtheorem{example}{Example}

\newtheorem{remark}{Remark}

\newtheorem{conjecture}{Conjecture}

\newtheorem{problem}{Problem}

\section{Some Convex Subsets of the Birkhoff Polytope}

A doubly stochastic (or bistochastic) matrix is a square non-negative matrix each of whose rows and whose columns sum to one.  We let $\mathcal{B}_n$ denote the set of $n$ by $n$ doubly stochastic matrices, $\mathcal{B}_n$ is sometimes called the Birkhoff polytope.  We will consider certain subsets of doubly stochastic matrices. The set of $n$ by $n$ symmetric bistochastic matrices, denoted $\mathcal{B}_n^{sym}$ is the set of all symmetric matrices in $\mathcal{B}_n$.  If $U$ is a $n$ by $n$ unitary matrix, let $U\circ \bar{U}$ be the $n$ by $n$ matrix whose $(i,j)$th entry is $|u_{ij}|^2$. $U\circ \bar{U}$ is clearly doubly stochastic; any doubly stochastic matrix that can be obtained in this way is called unistochastic.  (We follow standard notation is letting $A\circ B$ denote the Schur product of two $n$ by $n$ matrices $A$ and $B$; recall that $A\circ B$ is the $n$ by $n$ matrix whose $(i,j)$th entry is $a_{ij}b_{ij}$).  If $U$ is a Hermitian unitary matrix, then $U\circ \bar{U}$ is a symmetric bistochastic matrix.  Any matrix of the form $U\circ \bar{U}$ where $U$ is a Hermitian unitary matrix is called an $H$-unistochastic matrix.  $H$-unistochastic matrices were defined and first studied in \cite{st}.  When $n\geq 3$, there exist doubly stochastic matrices which are not unistochastic and symmetric bistochastic matrices which are not $H$-unistochastic.

\begin{example} \label{ex1} The matrix $ \left( \begin{smallmatrix} 0&\frac{1}{2}&\frac{1}{2} \\ \frac{1}{2}&0&\frac{1}{2} \\ \frac{1}{2}&\frac{1}{2}&0 \end{smallmatrix} \right)$ is an example of a matrix which is symmetric bistochastic but not unistochastic. \end{example}

 $\mathcal{B}_n$ and $\mathcal{B}_n^{sym}$ are both convex sets.  It is a well known result of Birkhoff that the extreme points of $\mathcal{B}_n$ are the permutation matrices.  Katz \cite{katz} proved that all extreme points of $\mathcal{B}_{n}^{sym}$ are of the form $\frac{1}{2}(P+P^T)$ where $P$ is a permutation matrix whose corresponding permutation contains no even cycle of length greater than or equal to four.  We state Katz's theorem more exactly below.

\begin{theorem} \cite{katz} \label{kthm} A matrix $M$ is an extreme point of the convex set of $n$ by $n$ symmetric bistochastic matrices if and only if there exists a permutation matrix $P$ such that $PMP^T$ is a block diagonal matrix with diagonal blocks $\{B_{n_k}\}_{k=1}^{m}$ where $m\geq 1$, $n_k$ is either odd or equal to two for each $k$, $\sum_{k=1}^{m}n_k=n$ and where $B_{n_k}$ is the $n_k$ by $n_k$ matrices defined as follows: $B_1=\left( 1\right)$, $B_2= \left( \begin{smallmatrix} 0& 1 \\ 1&0 \end{smallmatrix} \right)$ and for $n_k\ge 3$, $B_{n_k}$ is the matrix whose $(i,j)$ entry is $\frac{1}{2}$ if $i-j=\pm 1$ $(\mod n_k)$ and is zero otherwise.    \end{theorem}

We note that $B_3$ is the matrix given in example \ref{ex1}.  In fact for any odd $n$ greater than or equal to three, $B_n$ is a symmetric bistochastic matrix which is not unistochastic.  We will call the decomposition of an extreme point of $\mathcal{B}_n^{sym}$ given in the above theorem the Katz decomposition.

 Since every permutation matrix is clearly unistochastic, the set of $n$ by $n$ doubly stochastic matrices is the convex hull of the set of $n$ by $n$ unistochastic matrices.  When $n\geq 3$, the convex hull of the $n$ by $n$ the convex hull of the $H$-unistochastic matrices is strictly smaller than the set of symmetric bistochastic matrices, since any extreme point of $\mathcal{B}_n^{sym}$ whose Katz decomposition has blocks of size three by three or larger will not be unistochastic.

 We will let $\mathcal{B}_n^h$ denote the convex hull of the $n$ by $n$ $H$-unistochastic matrices.  We note that all three convex sets $\mathcal{B}_{n}$, $\mathcal{B}_n^{sym}$ and $\mathcal{B}_n^h$ are all invariant under permutation similarities (i.e. maps of the form $X\to PXP^T$).  The centroid of each of these convex sets is clearly the $n$ by $n$ matrix all of whose entries are $\frac{1}{n}$ as it is the unique point in $\mathcal{B}_n$ which is fixed by all permutation similarities. We will follow \cite{st} and let $W_n$ denote the $n$ by $n$ matrix all of whose entries are $\frac{1}{n}$.  We will examine the relationship between the convex sets $\mathcal{B}_n^{sym}$ and $\mathcal{B}_n^h$. We then apply the theory of completely positive linear maps to some of our results; we review the basics of completely positive linear maps in the next section.

\section{Introduction to Completely Positive Linear Maps}

In this section we review the basic theory of completely positive linear maps.

\begin{definition}\cite{paulsen}  A linear map $\Phi:\mathcal{M}_n(\mathbb{C})\rightarrow\mathcal{M}_n(\mathbb{C})$ is
\begin{enumerate}
\item positive if $\Phi(A)$ is positive semidefinite whenever $A$ is positive semidefinite.
\item completely positive if $\Phi\otimes I_{k}$ is a positive linear map from $\mathcal{M}_{kn}(\mathbb{C}) \to \mathcal{M}_{kn}(\mathbb{C})$ for all $k\in\mathbb{N}$.
\item trace-preserving if $Tr\left(\Phi(A)\right)=Tr(A)$, for all $A\in \mathcal{M}_{n}$.
\item unital if $\Phi(I_{n})=I_{n}$.
\end{enumerate}

We will call a unital completely positive trace preserving linear map a doubly stochastic map.

\end{definition}

The following theorem of Choi provides a necessary and sufficient condition on the decomposition of linear map to be completely positive.

\begin{theorem}\cite{choi}
A linear map $\Phi:\mathcal{M}_n(\mathbb{C})\longrightarrow\mathcal{M}_n(\mathbb{C})$ is completely positive if and only if it admits the following form
\begin{eqnarray}
\Phi(A)=\sum_{i=1}^{r}V_iAV_i^*
\end{eqnarray}
where $\left\{V_i\right\}\in\mathcal{M}_{n}(\mathbb{C})$.
\end{theorem}

The matrices $\left\{V_i\right\}$ are called the Kraus operators for a given map $\Phi$, and the choice of such Kraus operators of $\Phi$ is non-unique.  The number of terms in the summation will be at a minimum whenever the Kraus operators are chosen to be linearly independent, this minimal value of $r\leq n^2$ is called the rank of $\Phi$.  The rank one doubly stochastic maps are precisely the unitary transforms.

\begin{definition}
  Let $\Phi:\mathcal{M}_n(\mathbb{C})\longrightarrow\mathcal{M}_n(\mathbb{C})$ be a completely positive map. $\Phi$ is a mixed unitary map if there exists a collection of unitary matrices $\left\{U_i\right\}_{i=1}^m$ such that
  \begin{eqnarray*}
    \Phi(A)=\sum_{i=1}^mp_iU_iAU_i^*
  \end{eqnarray*}
  where $\sum_{i=1}^mp_i=1$ and $p_i>0$.
\end{definition}

Mixed unitary maps are also often called random unitary maps in the literature; in using the term mixed unitary we are following \cite{watrous}.
Landau and Streater \cite{streater} have shown that every doubly stochastic map on $\mathcal{M}_2(\mathbb{C})$ is a mixed unitary map but there exist doubly stochastic maps on $\mathcal{M}_n(\mathbb{C})$ for all $n\ge 3$ which are not mixed unitary.  Much more about doubly stochastic and mixed unitary maps and their application to physics can be found in \cite{wolf}.

\section{Geometrical Relations between some convex subsets of the Birkhoff polytope}

The following general problem has been considered in many different contexts.

\begin{problem} \label{centr} Let $V$ be a finite dimensional real of complex vector space and let $A$ and $B$ be two compact convex subsets of $V$ with $B\subseteq A$. Suppose $A$ and $B$ have a common centroid $c\in B$.  Find $\lambda$, the largest constant in $[0,1]$ for which the set inclusion $\lambda A +(1-\lambda)c \subseteq B$ is satisfied. \end{problem}

Some norm inequalities are in fact essentially the solution to this problem where $A$ and $B$ are the unit balls of the norms in question.  This problem has also been explored in the context of doubly stochastic maps and correlation matrices.
 In \cite{watrous}, the case where $A$ is the set of doubly stochastic maps on $M_{n}(\mathbb{C})$ and $B$ is the set of mixed unitary maps on $M_{n}(\mathbb{C})$ was considered.  Here a lower bound of $\frac{1}{n^2-1}$ was found for $\lambda$ for each fixed $n$.    In \cite{marshall}, the case where $A$ is the $n$ by $n$ real correlation matrices and where $B$ is the the convex hull of the $n$ by $n$ real rank one correlation matrices was considered.  It was shown that $\lambda=\frac{2}{3}$ when $3\le n\le 6$ and bounds were given for higher dimensional cases.  The motivations for the latter two problems are quantum information and semidefinite programming respectively.

We will consider the case of problem \ref{centr} where $A$ is the set of $n$ by $n$ symmetric bistochastic matrices and $B$ is the convex hull of the $n$ by $n$ $H$-unistochastic matrices.

\begin{theorem} Let $\lambda_n$ be the largest constant in $[0,1]$ for which $\lambda_n M+(1-\lambda_n)W_n\in \mathcal{B}_n^h$ for all $n$ by $n$ symmetric bistochastic matrices $M$.  Then $\lambda_1=\lambda_2=1$, $\lambda_3=\lambda_4=\frac{2}{3}$. \end{theorem}

That $\lambda_1=\lambda_2=1$ is an immediate consequence of the fact that all $1$ by $1$ and $2$ by $2$ symmetric bistochastic matrices are $H$-unistochastic. We split the remainder of this result into two parts as the proofs are quite different for each case:

\begin{proposition} $\lambda_3=\frac{2}{3}$.  \end{proposition}

\noindent\emph{Proof:}  We note that there is only one extreme point of $\mathcal{B}_{3}^{sym}$ which is not $H$-unistochastic; this is the matrix $M=\left( \begin{smallmatrix} 0&\frac{1}{2}&\frac{1}{2} \\ \frac{1}{2}&0&\frac{1}{2} \\ \frac{1}{2}&\frac{1}{2}&0 \end{smallmatrix} \right) $ from example \ref{ex1}.  So we only need show that $kM+(1-k)W_3$ is in $H_3$ for $k=\frac{2}{3}$ but not in $H_3$ if $k>\frac{2}{3}$.  Note that $\frac{2}{3}M+\frac{1}{3}W_3=\left( \begin{smallmatrix} \frac{1}{9}&\frac{4}{9}&\frac{4}{9} \\ \frac{4}{9}&\frac{1}{9}&\frac{4}{9} \\ \frac{4}{9}&\frac{4}{9}&\frac{1}{9} \end{smallmatrix} \right) $ is $H$-unistochastic since $\left( \begin{smallmatrix} -\frac{1}{3}&\frac{2}{3}&\frac{2}{3} \\ \frac{2}{3}&-\frac{1}{3}&\frac{2}{3} \\ \frac{2}{3}&\frac{2}{3}&-\frac{1}{3} \end{smallmatrix} \right) $ is a Hermitian unitary matrix.

We note that if $n$ is odd and $U$ is an $n$ by $n$ Hermitian unitary matrix, then $\vert tr(U)\vert \ge 1$ since all the eigenvalues of $U$ are plus or minus one.  Hence by the Cauchy-Schwarz inequality, $\sum_{i=1}^{n}\vert u_{ii} \vert^2\ge \frac{1}{n}\vert \sum_{i=1}^{n} u_{ii} \vert^2 \ge \frac{1}{n}$ which implies that the trace of any matrix in $\mathcal{B}_n^h$ must be greater than or equal to $\frac{1}{n}$ when $n$ is odd.  Since $tr(kM+(1-k)W_3)=1-k\le \frac{1}{3}$ when $k\ge \frac{2}{3}$ our result follows. $\blacksquare $

\begin{proposition} $\lambda_4= \frac{2}{3}$.  \end{proposition}

\noindent \emph{Proof:} We note that up to permutation similarity there is only one extreme point of $\mathcal{B}_{4}^{sym}$ which is not $H$-unistochastic; this is the matrix $M=\left( \begin{smallmatrix} 0&\frac{1}{2}&\frac{1}{2}&0 \\ \frac{1}{2}&0&\frac{1}{2}&0 \\ \frac{1}{2}&\frac{1}{2}&0&0 \\ 0&0&0&1 \end{smallmatrix} \right) $.  Now let $X=\left( \begin{smallmatrix} 0&\frac{1}{3}&\frac{1}{3}&\frac{1}{3} \\ \frac{1}{3}&0&\frac{1}{3}&\frac{1}{3} \\ \frac{1}{3}&\frac{1}{3}&0&\frac{1}{3} \\ \frac{1}{3}&\frac{1}{3}&\frac{1}{3}&0 \end{smallmatrix} \right) $. This matrix is arithmetic mean of the three Hermitian permutation matrices corresponding to the three permutations $(12)(34)$, $(13)(24)$ and $(14)(23)$ and hence $X\in \mathcal{B}_4^h$. Let $Y=\left( \begin{smallmatrix} \frac{1}{9}&\frac{4}{9}&\frac{4}{9}&0 \\ \frac{4}{9}&\frac{1}{9}&\frac{4}{9}&0 \\ \frac{4}{9}&\frac{4}{9}&\frac{1}{9}&0 \\ 0&0&0&1 \end{smallmatrix} \right) $; it follows from the proof of the previous proposition that $Y\in \mathcal{B}_4^h$.  We note that $\frac{2}{3}M+\frac{1}{3}W_4=\frac{1}{4}X+\frac{3}{4}Y\in \mathcal{B}_4^h$.  Hence $\lambda_4\ge \frac{2}{3}$.

We now prove that if $A\in \mathcal{B}_4^h$, then $3a_{11}+3a_{22}+3a_{33}-a_{44}\geq 0$. We only need show this for the special case where $A$ is $H$-unistochastic; the general case where $A$ is in the convex hull of $H$-unistochastic matrices follows from convexity. Let $U$ be a Hermitian unitary matrix such that $U\circ \bar{U}=A$.  Note that $tr(U)$ is an even integer and $u_{44}\in [-1,1]$, hence $\vert u_{44}\vert \le \vert u_{44}-tr(U)\vert = \vert u_{11}+u_{22}+u_{33}\vert \le \sum_{i=1}^{3}\vert u_{ii}\vert$.  Squaring both sides of the last inequality and then using Cauchy-Schwarz, we get $a_{44}=\vert u_{44}\vert^2 \le (\sum_{i=1}^{3}\vert u_{ii}\vert)^2 \le 3\sum_{i=1}^{3}\vert u_{ii}\vert^2=3\sum_{i=1}^{3}a_{ii}$. Now the diagonal entries of $kM+(1-k)W_4$ are $(\frac{1-k}{4},\frac{1-k}{4},\frac{1-k}{4},\frac{1+3k}{4})$.  Hence if $kM+(1-k)W_4\in \mathcal{B}_{4}^{h}$, then $9-9k\ge 1+3k$ which means that $k\le \frac{2}{3}$. Therefore $\lambda_4=\frac{2}{3}$. $\blacksquare $

We note that the number $\frac{2}{3}$ appears both as a solution to the small dimensional cases of both the correlation matrix problem from \cite{marshall} and our problem on symmetric bistochastic matrices.  While this may be coincidental, we will soon derive a result which suggests some connections between these two problems.  We first need to introduce the class of self-dual completely positive linear maps and explore its connections to both the symmetric bistochastic matrices and the real correlation matrices.  We do this in the next section, which culminates in a partial result of our problem for general $n$.  The exact value of $\lambda_n$ for $n \geq 5$ remains an open problem.

\section{Self-dual linear maps and correlation matrices}

As the names suggest, there are very close connections between doubly stochastic matrices and doubly stochastic maps.  (Some of these connections will be reviewed later on in the paper).  The doubly stochastic maps are in a sense an analogue of the doubly stochastic matrices.  We would like to find a similar analogue of the symmetric bistochastic matrices.
We note that $\mathcal{M}_n(\mathbb{C})$ is a Hilbert space under the inner product $\langle A,B \rangle= Tr(AB^*)$.  The adjoint or dual to the CP map $\Phi$ is the linear map $\Phi^*$ which satisfies $\langle A, \Phi(B) \rangle = \langle \Phi^*(A), B\rangle$ for all $A,B$.  It follows from simple calculation that the dual of $\Phi(A)=\sum_{i=1}^{r}V_iAV_i^*$ is  $\Phi^*(A)=\sum_{i=1}^{r}V_i^*AV_i$. Hence the dual of a CP map is a CP map.  We say that a CP map is self-dual if $\Phi=\Phi^*$.  The self-dual doubly stochastic maps are the analogue of the symmetric bistochastic matrices.  We study the properties of self-dual completely positive linear maps in this section.

The following theorem is straightforward and can be found in \cite{mwolftext}.

\begin{theorem}\label{selfdualhermkraus} \cite{mwolftext}
A completely positive map $\Phi:\mathcal{M}_n(\mathbb{C})\longrightarrow\mathcal{M}_n(\mathbb{C})$ is self-dual if and only if there exists a set of $n$ by $n$ Hermitian matrices $\left\{H_i\right\}$ such that:

\begin{eqnarray}
\Phi(A)=\sum_{i=1}^{r}H_iAH_i
\end{eqnarray}

for all $A\in\mathcal{M}_n(\mathbb{C})$.
\end{theorem}

\noindent\emph{Proof:}
The if direction is immediate. Now suppose $\Phi(A)=\sum_{i=1}^{r}V_iAV_i^*$ is self-dual.
Since $\Phi=\Phi^*$, the self-dual map $\Phi$ can be expressed as $\Phi=\frac{\Phi+\Phi^*}{2}$. Since any matrix V can be decomposed via Cartesian decomposition as $V=K+iL$, where K and L are $n\times n$ Hermitian matrices, decomposing each Kraus operator in this manner yields the result $\Phi(A)=\sum_{i=1}^{r}K_iAK_i+L_iAL_i$ $\blacksquare$\\

We recall the result that two sets of Kraus operators $\{V_i\}$ and $\{W_j\}$ represent the same completely positive map $\Phi$ if and only if there exists a unitary matrix $U=(u_{ij})$ such that $V_i=\sum_j u_{ij}W_j$ \cite{choi}.
While not every set of Kraus operators of a self-dual map will be Hermitian, theorem \ref{selfdualhermkraus} guarantees that there is always at least one set of Kraus operators which are Hermitian for a self-dual map.  We now introduce a self-dual version of the concept of a mixed unitary map.

\begin{definition}
  Let $\Phi:\mathcal{M}_n(\mathbb{C})\longrightarrow\mathcal{M}_n(\mathbb{C})$ be a self-dual completely positive map. $\Phi$ is a mixed Hermitian unitary map if there exists a collection of Hermitian unitary matrices $\left\{U_i\right\}_{i=1}^m$ such that
  \begin{eqnarray*}
    \Phi(A)=\sum_{i=1}^mp_iU_iAU_i
  \end{eqnarray*}
  where $\sum_{i=1}^mp_i=1$ and $p_i>0$.
\end{definition}

 The relationship between the mixed Hermitian unitary maps and the self-dual doubly stochastic maps can be thought of as being analogous to the relationship between the mixed unitary maps and the doubly stochastic maps. There is one key difference which we will describe in the next paragraph.

 Let $\Phi$ be a doubly stochastic map on $\mathcal{M}_n(\mathbb{C})$, and define $\Delta_{\Phi}(D)=I_{n}\circ \Phi(D)$ for all $n$ by $n$ diagonal matrices $D$.  Since the $n$ by $n$ diagonal matrices form an $n$ dimensional vector space with the obvious standard basis, we will identify $\Delta_{\Phi}$ with its $n$ by $n$ matrix representation.  Hence if $\Phi$ is a doubly stochastic map, then $\Delta_{\Phi}$ is a doubly stochastic matrix.  If $\Phi$ is a self-dual doubly stochastic map, then $\Delta_{\Phi}$ is a symmetric bistochastic matrix.   Conversely, if $M$ is a doubly stochastic matrix (resp. symmetric bistochastic matrix), there exists a doubly stochastic map (resp. self-dual doubly stochastic map) $\Phi$ such that $\Delta_{\Phi}=M$.  The converse statements follows from considering the rank one map $\Phi(A)=PAP^T$ where $P$ is any permutation matrix. Then it can be shown that $\Delta_{\Phi}=P$. By linearity, maps of rank $r\leq n^2$ must be doubly stochastic. If $U$ is a unitary matrix and $\Phi(A)=UAU^*$, then $\Delta_{\Phi}=U \circ \bar{U}$; hence by linearity of $\Delta$ we have the following result:

\begin{proposition} \label{useful} Let $\Phi$ be a completely positive linear map on $\mathcal{M}_n(\mathbb{C})$.  If $\Phi$ is mixed Hermitian unitary, then $\Delta_{\Phi}$ lies in the convex hull of the $n$ by $n$ $H$-unistochastic matrices.  \end{proposition}

 This gives us a key difference between the self-dual case and the general case.  Given a doubly stochastic linear map $\Phi$, the matrix $\Delta_{\Phi}$ will never give us any information as to whether $\Phi$ is mixed unitary since the convex hull of the $n$ by $n$ unistochastic matrices is the entire set of $n$ by $n$ doubly stochastic matrices.

We introduce the set of correlation matrices which will be useful in constructing doubly stochastic maps.

\begin{definition}
  Let $C\in\mathcal{M}_n(\mathbb{C})$ be a positive semidefinite matrix with $c_{ii}=1$ for all $1\leq i\leq n$. Then we say that C is a correlation matrix.
\end{definition}

 For any $n$ by $n$ correlation matrix $C$, we follow the convention in \cite{LW} by letting $\Phi_C$ denote the completely positive map which maps $A\in M_n$ to the Schur product $A\circ C$.  We note that if $C$ is an $n$ by $n$ correlation matrix, then $\Phi_C$ is a doubly stochastic map with $\Phi_{C^T}$ being its dual. Hence $\Phi(A)=A\circ C$ is self-dual if and only if $C$ is a real correlation matrix.  Any completely positive linear map which leaves every diagonal matrix unchanged must be of this type; this statement of the result was given independently in \cite{kye} and \cite{LW}. See for instance proposition 2.1 of Kye \cite{kye} and the discussion previous to it or theorem 1 from \cite{LW}.  A more general result was found even earlier in \cite{uhlmann} and was rediscovered in \cite{Chefles200014}. We will not use the more general result here and refer to \cite{Chefles200014, cjw, uhlmann} for details.

\begin{proposition}\label{kyeresult} \cite{kye,LW}
  Let $\Phi$ be a completely positive linear map on the $n$ by $n$ matrices such that $\Phi(D)=D$ for all diagonal matrices $D$. Then there exists an $n$ by $n$ correlation matrix C such that $\Phi(A)=\Phi_C(A)=C\circ A$ for all $A\in\mathcal{M}_n(\mathbb{C})$.
\end{proposition}

The next two results provide conditions on whether such maps are mixed unitary or mixed Hermitian unitary.

\begin{theorem} \label{b1}
Let $C$ be a correlation matrix. Then $\Phi_C$ is a mixed unitary map if and only if the associated correlation matrix $C$ is a convex combination of rank 1 correlation matrices.
\end{theorem}

\noindent \emph{Proof:} Suppose $\Phi(A)$ is a mixed unitary map, then $\Phi(A)=\sum_{k=1}^rp_kU_kAU_k^*$ where $\{ U_k\}_{k=1}^{r}$ are unitary, $\sum_{k=1}^r p_k=1$ and $p_k>0$.  Let $E_{i}$ be the matrix with a one in its $i$th row, $i$th column entry and zeros everywhere else.  Since $\Phi(E_{i})=E_{i}$ which is rank one, $U_kE_{i}U_k^*=E_{i}$ for $1\le i \le n$ and $1\le k \le n$.  Therefore $U_k$ commutes with all diagonal matrices which means it must be a diagonal matrix. Let $v_k$ be the vector consisting of the diagonal elements of $U_k$, then $v_kv_k^{*}$ is a rank-one correlation matrix and $C=\sum_{k=1}^{r}p_kv_kv_k^*$. Similarly if $C=\sum_{k=1}^{r}p_kv_kv_k^*$ where $v_kv_k^{*}$ are rank-one correlation matrices,  then we let $U_k$ be the diagonal unitary matrices formed using the entries of $v_k$ and hence $C\circ A=\sum_{k=1}^rp_kU_kAU_k^*$. $\blacksquare$

By using the same argument as above, we obtain:

 \begin{theorem}\label{b2}
Let $C$ be a real correlation matrix. Then $\Phi_C$ is a mixed Hermitian unitary map if and only if the associated correlation matrix $C$ is a convex combination of real rank 1 correlation matrices.
\end{theorem}

In light of this, it would be useful to be able to test whether a real correlation matrix $C$ is the convex hull of the real rank one correlation matrices.  The following result is a very slight modification of a result from \cite{guptpar}.

\begin{theorem} \label{BGP} Let $C$ be an $n$ by $n$ real correlation matrix and let $P_n$ be the collection of all subsets of $\{ 1,2,...,n\}$.  Then $C$ is in the convex hull of the real rank one correlation matrices if and only if there exists a function $f:P_n\to \mathbb{R}$ with the following properties:
\begin{enumerate}
\item $f(\emptyset )=1$
\item $f(\{i,j\} )=c_{ij}$ for $1\le i<j \le n$
\item  $\sum_{T\in P_n} (-1)^{\vert S\cap T \vert} f(T)\geq 0$ for all $S\in P_n$
\end{enumerate} \end{theorem}

\noindent \emph{Proof:} Let $K$ be the set of all functions $f:P_n\to \mathbb{R}$ which satisfy conditions one and three from the statement of the theorem above. It is clear that $K$ is a convex set.  For all $A\in P_n$, let $f_A (T)=(-1)^{\vert A\cap T \vert}$.  Since $\sum_{T\in P_n} (-1)^{\vert S\cap T \vert} f_A(T)=0$ except when $S=A$, $\{f_A \}_{A\in P_n}$ are exactly the extreme points of $K$.  Let $C_A$ be the $n$ by $n$ matrix whose main diagonal elements are all one and whose $(i,j)$th entry is $f_A (\{ i,j\} )$ for all $i\neq j$.  Since $\{C_A \}_{A\in P_n}$ is the set of all $n$ by $n$ real rank one correlation matrices, the result follows. $\blacksquare $

 We can use Theorem \ref{BGP} to give us a lower bound for a rank-dependent variant of the problem from \cite{marshall}.

 \begin{corollary}\label{coro} Let $C$ be an $n$ by $n$ real correlation matrix and let $r=rank(C)$.  Then $\frac{1}{r}C+\frac{r-1}{r}I$ is in the convex hull of real rank one correlation matrices.  \end{corollary}

 \noindent \emph{Proof:} There exists a real $r$ by $n$ matrix $A$ such that $C=A^T A$. We note that all entries of $A$ must have absolute value less than or equal to one. Now we construct the following function $g:P_n\to \mathbb{R}$.  Let $g(\emptyset )=r$ and let $g(T)=\sum_{i=1}^{r} \prod_{j\in T}a_{ij}$ whenever $T$ is non-empty.  We note that $g(\{ i,j\} )=c_{ij}$ and for any $S\in P_n$, $\sum_{T\in P_n} (-1)^{\vert S\cap T \vert} g(T)=\sum_{i=1}^{r}\sum_{T\in P_n} (-1)^{\vert S\cap T \vert} \prod_{j\in T}a_{ij}$ $=\sum_{i=1}^{r}\prod_{j=1}^{n}(1+(-1)^{\vert \{j\} \cap S\vert }a_{ij})\ge 0$.  Now take the function $f(T)=r^{-1}g(T)$ and apply $f$ to theorem \ref{BGP} to show that $\frac{1}{r}C+\frac{r-1}{r}I$ is in the convex hull of real rank one correlation matrices.

The set of all $n$ by $n$ correlation matrices in $\mathcal{M}_n(\mathbb{F})$ over a field $\mathbb{F}$ which is either the real or the complex numbers form a compact convex set that is usually denoted $\varepsilon_{n}(\mathbb{F})$.  As with most convex sets, there has been interest in the properties of its extreme points. It was shown by Christensen and Vesterstrom \cite{christ} and Loewy \cite{loewy} that  $\varepsilon_{n}(\mathbb{C})$ has extreme points of rank r if and only if $r\leq \sqrt n$. Later, Grone, Pierce and Watkins  \cite{grone} found an analogous result for the set $\varepsilon_{n}(\mathbb{R})$ which they showed has extreme points of rank $r$ if and only if $r(r+1)\leq 2n$.  A simple derivation of these results along with an efficient algorithm for checking whether a correlation matrix is extremal was given by Li and Tam in \cite{LiTam}.  The Li-Tam characterization of extremal correlation matrices has been generalized to the infinite-dimensional case by Kiukas and Pellonp\"{a}\"{a} in \cite{kp}.

 We let $\rho_{n,r}$ be the smallest number positive number less or equal to one than one for which $\rho_{n,r}C+(1-\rho_{n,r})I_n$ is in the convex hull of rank one correlation matrices whenever $C$ is an $n$ by $n$ real correlation matrix of rank at most $r$.  From \cite{marshall}, we know that $\rho_{3,2}=\rho_{4,2}=\rho_{5,2}=\rho_{6,2}=\rho_{6,3}=\frac{2}{3}$.  We also know that when $n\geq 3$ and $k\geq 2$,  $\frac{1}{r} \le \rho_{n,r} \le \frac{2}{3}$ where the lower bound follows from corollary \ref{coro} and the upper bound follows from the fact that $\{ \rho_{n,r} \}$ is a nonincreasing sequence in both $n$ and $r$.

We now return to our main problem.  First we recall some notation.  For any natural number $n$, let $F_n$ denote the $n$ by $n$ Fourier matrix which is the unitary matrix whose $(j,k)$th entry is $n^{-\frac{1}{2}}\e^{\frac{2\pi (j-1)(k-1)}{n}i}$.  If $U$ is any $n$ by $n$ unitary, then $\Gamma_{U}$ denotes the unitary conjugation map which sends $A\in M_n$ to $UAU^*$.

\begin{proposition} \label{part} Let $n\geq 3$, let $m,q\in \mathbb{N}$ with $mq=n$ and let $M$ be a $n$ by $n$ be an extremal symmetric bistochastic matrix having $q$ Katz blocks all of which are $m$ by $m$. Then $\rho_{n,2} M+(1-\rho_{n,2}) W_n\in \mathcal{B}_n^h$. \end{proposition}

\noindent\emph{Proof:} In this proof we note that we are using the symbol $\circ$ to denote composition and not Schur multiplication.  Let $n=mq$ and let $M$ be the extremal symmetric bistochastic matrix having $q$ Katz blocks all of which are $m$ by $m$. Let $\omega=e^{\frac{2\pi i}{m}}$ and $q=\frac{n}{m}$.  Then $M$ is permutationally similar to $B_m\otimes I_q$, so without loss of generality let us assume $M=B_m\otimes I_q$.

Let $C_{m}$ be the $m$ by $m$ correlation matrix whose $(j,k)$th entry is $\cos (\frac{2\pi(j-k)}{m})=\frac{\omega^{j-k}+\omega^{k-j}}{2}$.  We note that the $(j,k)$th entry of $C_{m}$ is equal to $\langle e_j, e_k \rangle$ where $e_k=(\cos (\frac{2\pi k}{m}),\sin (\frac{2\pi k}{m}))^T$ is the unit vector with polar angle $\frac{2\pi k}{m}$ radians for all $k$; hence $C_m$ is a rank two correlation matrix.  Let $v_p$ be the $p$th column of $F_m$, then $v_pv_p^*$ is an $m$ by $m$ matrix whose $(j,k)$th entry is $\frac{1}{m}\omega^{(p-1)(j-k)}$. Hence $\Phi_{C_m}(v_pv_p^*)=\frac{1}{2}(v_{p-1}v_{p-1}^*+v_{p+1}v_{p+1}^*)$ where the additions $p-1$ and $p+1$ are taken mod $p$. Let $E_{p}$ denote the $m$ by $m$ matrix with a one in its $(p,p)$th entry and zeros everywhere else. Note that $v_{p}v_{p}^*=F_mE_pF_m^*$.  Let $\Phi=\Gamma_{F_{m}^*}\circ \Phi_{C_m} \circ \Gamma_{F_{p}}$ , then $\Phi_n(E_p)=\frac{1}{2}(E_{p-1}+E_{p+1})$ where the additions $p-1$ and $p+1$ are taken mod $m$.

Note that $\Delta_{\Phi_m\otimes id_q}=M$.  Let $U=F_m\otimes F_q$, then $\Phi_m\otimes id_q=\Gamma_{U^*}\circ\Phi_{qC_m\otimes W_q} \circ \Gamma_U$.
Let $\Xi_n=\Gamma_{U^*}\circ\Phi_{I_n} \circ \Gamma_U$. We note that for any $n$ by $n$ diagonal matrix $D$,  $\Xi(D)=\frac{tr(D)}{n}I_n$ and so $\Delta_{\Xi_n}=W_n$.  Since $qC_m\otimes W_q$ is an $n$ by $n$ rank two correlation matrix, $B=\rho_{n,2}(qC_m\otimes W_q)+(1-\rho_{n,2}) I_n$ is a convex combination of rank one matrices.  Hence by theorem \ref{b2}, $\Psi=\rho_{n,2}(\Phi_m\otimes id_q)+(1-\rho_{n,2})\Xi_n=\Gamma_{U^*}\circ\Phi_{B} \circ \Gamma_U$ is mixed Hermitian unitary map. Therefore $\rho_{n,2} M+(1-\rho_{n,2}) W_n=\Delta_{\Psi}\in \mathcal{B}_n^h$ by proposition \ref{useful}.

\section{The self-dual Landau-Streater theorem}

The set of doubly stochastic maps on $\mathcal{M}_n(\mathbb{C})$ is convex and as with most convex sets, there has been interest in finding the extreme points.  Since the unitary transforms $\Gamma_{U}(A)=UAU^*$ are the rank one doubly stochastic maps, they must be extremal.  Are they the only extreme points?
Landau and Streater \cite{streater} have shown that every extreme doubly stochastic map on $\mathcal{M}_2(\mathbb{C})$ must be a unitary transform, but there exist extremal doubly stochastic maps on $\mathcal{M}_n(\mathbb{C})$ which are not unitary transforms for all $n$ greater than or equal to three.  It follows from this that every doubly stochastic map on $\mathcal{M}_2(\mathbb{C})$ is a mixed unitary map but there exist doubly stochastic maps on $\mathcal{M}_n(\mathbb{C})$ for all $n\ge 3$ which are not mixed unitary.  In this section, we show that the situation for self-dual maps is similar; every self-dual doubly stochastic map on $\mathcal{M}_2(\mathbb{C})$ is a mixed Hermitian unitary map but there exists doubly stochastic maps on $\mathcal{M}_n(\mathbb{C})$ for all $n\ge 3$ which are not mixed Hermitian unitary.  Our results on correlation matrices give us a particularly simple construction of the latter.

\begin{example} \label{exam2} Let $C= \left( \begin{smallmatrix} 1&-\frac{1}{2}&-\frac{1}{2} \\ -\frac{1}{2}&1&-\frac{1}{2} \\ -\frac{1}{2}&-\frac{1}{2}&1 \end{smallmatrix} \right) $.  $C$ is a real correlation matrix which is in the convex hull of complex rank-one correlation matrices but not in the convex hull of the real rank-one correlation matrices. (In fact, $C$ is a rank 2 extreme point of $\varepsilon_{3}(\mathbb{R})$). Hence the map $A\longrightarrow A\circ C$ is a self-dual mixed unitary map which is not mixed Hermitian unitary. \end{example}

Since the convex set of $n$ by $n$ real correlation matrices has rank two extreme points for all $n\geq 3$, there exist self-dual doubly stochastic maps on $\mathcal{M}_n(\mathbb{C})$ which are not mixed Hermitian unitary.  We can also use the fact that there are symmetric bistochastic matrices which are not the convex combination of $H$-unistochastic matrices to give another class of examples.

\begin{example} Let $n\geq 3$ and let $P$ be an $n$ by $n$ permatation matrix corresponing to the permutation $(123)$.  Let $\Phi(A)=\frac{1}{2}(PAP^T+P^TAP)$.  Then $\Delta_{\Phi}$ consists of the extremal symmetric bistochastic matrix which is the direct sum of a single Katz block $B_3$ with $n-3$ Katz blocks $B_1$.
Since $\Delta_{\Phi}$ is not in the convex hull of the $H$-unistochastic matrices, $\Phi$ cannot be mixed Hermitian unitary.  Note however that this map like that of example \ref{exam2} is both self-dual and mixed unitary.
    \end{example}

 We now prove that for $n=2$, all the self-dual doubly stochastic maps are mixed Hermitian unitary.  This result is the self-dual analogue of the Landau-Streater theorem.

\begin{theorem}\label{randomrefthm}
Let $\Phi:\mathcal{M}_2(\mathbb{C})\longrightarrow\mathcal{M}_2(\mathbb{C})$ be a doubly stochastic map. Then $\Phi$ is self-dual if and only if it admits a Kraus decomposition of the form
\[\Phi(A)=\sum_i^m p_iR_iAR_i\]
where $\{R_i\}_{i=1}^m$ is a collection of unitary Hermitian matrices and $\sum_{i=1}^{m}p_i=1$ with $p_i>0$.
\end{theorem}


\noindent\emph{Proof:} The if direction is immediate.
We prove the only if part.
Let $\Psi$ be an extreme doubly stochastic map on $\mathcal{M}_2(\mathbb{C})$. Then we will show that the map $\Phi=\tfrac{\Psi+\Psi^*}{2}$ admits a Kraus decomposition composed of Hermitian unitary matrices and the general case will follow by linearity. Since the extreme points are exactly the unitary transforms, let
$U=VDV^*$ be the spectral decomposition of the single Kraus operator $U$ such that $\Phi(A)=\phi\circ\Psi_D\circ\phi^*(A)$,
where $\phi:A\longrightarrow VAV^*$ and $\Psi_D:A\longrightarrow \frac{DAD^*+D^*AD}{2}$. Clearly, the diagonal matrix resulting from the spectral decomposition is $D=\mbox{diag}(e^{i\theta},e^{-i\phi})$. By Proposition \ref{kyeresult}, $\Psi_D$ is a Schur map expressed as
\begin{eqnarray*}
\Psi_D(A)&=&
\left[\begin{array}{cc}
1 & \cos{(\theta+\phi)} \\
\cos{(\theta+\phi)} & 1 \\
\end{array}\right]\circ A\\
\end{eqnarray*}
 Since every $2\times 2$ real correlation matrix is a convex combination of real rank 1 matrices, it is clear that $\Phi(A)=\phi\circ\Psi\circ\phi^*(A)=\sum_{i=1}^qp_i(VD_iV^*)A(VD_iV^*)=\sum_{i=1}^qp_iR_iAR_i$ where $\mbox{diag}(D_i)=v_i$ such that $W_i=v_iv_i^T$. $\blacksquare$\\

\section{Conclusion $\&$ Outlook}

We conclude this paper by summarizing our main problem and its possible connection with the correlation matrix problem studied in \cite{marshall}.  We summarize the relationships between the mathematical objects studied in the paper and formulate a conjecture based on our partial results.

While the convex hull of the $n$ by $n$ unistochastic matrices is the set of $n$ by $n$ doubly stochastic matrices, $\mathcal{B}_n^{h}$ the convex hull of the $n$ by $n$ $H$-unistochastic matrices is a proper subset of $\mathcal{B}_n^{sym}$ the set of symmetric bistochastic matrices.  We considered the problem of finding the smallest nonnegative number $\lambda_n$ for which the set $\lambda_n \mathcal{B}_n^{sym} + (1-\lambda_n) W_n$ is contained in the set $\mathcal{B}_n^{sym}$.  Clearly $\lambda_1=\lambda_2=1$.  We found the first two non-trivial cases of $\lambda_n$ by showing that $\lambda_3=\lambda_4=\frac{2}{3}$.

In \cite{marshall}, a similar problem was considered. Let $\varepsilon_{n}(\mathbb{R})$ be the convex set of all $n$ by $n$ real correlation matrices and let $\mathcal{Q}_n$ be the convex hull of the set of $n$ by $n$ real rank one correlation matrices.  The problem in \cite{marshall} is to find the smallest nonnegative number $\rho_n$ for which the set $\rho_n \varepsilon_{n}(\mathbb{R})+(1-\rho_n)I_n$ is completely contained in $\mathcal{Q}_n$.  The known values of $\rho_n$ are $\rho_1=\rho_2=1$ and $\rho_3=\rho_4=\rho_5=\rho_6=\frac{2}{3}$\cite{marshall}.  We note that $\lambda_n=\rho_n$ when $n\leq 4$ which accounts for all of the known values of $\lambda_n$.  There are also other connections between these two problems.    Specifically, $\Phi_C(A)=C\circ A$ is a self-dual doubly stochastic map if and only if $C$ is a real correlation matrix and $\Phi_C$ is mixed Hermitian unitary if and only if $C$ is in $\mathcal{Q}_n$.  Furthermore, if $\Phi$ is any self-dual doubly stochastic map, then $\Delta_{\Phi}$ is a symmetric stochastic matrix; if in addition $\Phi$ is mixed Hermitian unitary then $\Delta_{\Phi}\in\mathcal{B}_n^{h}$. These two results were used to prove proposition \ref{part} which also suggests a connection between $\lambda_n$ and $\rho_n$.  Motivated by this, we offer the following conjecture:

\begin{conjecture} $\lambda_n=\rho_n$ for all $n\in \mathbb{N}$. \end{conjecture}

In the previous section, as a further application of these concepts, we show that every self-dual doubly stochastic maps on $M_{n}(\mathbb{C})$ is mixed Hermitian unitary if and only if $n\leq 2$.  This is the self-dual version of the well-known result of Landau and Streater that every doubly stochastic maps on $M_{n}(\mathbb{C})$ is mixed unitary if and only if $n\leq 2$.

\section{Acknowledgements}

R.P. was supported by the NSERC Discovery grant number 400550 and would like to thank Murray Marshall for several discussions on correlation matrices.  Both authors would like to thank the anonymous referee and Prof. Armin Uhlmann for suggesting many useful improvements to this paper.

\bibliography{thesisn}
\bibliographystyle{plain}

\end{document}